\documentclass[11pt]{article}
\usepackage{mathrsfs}
\usepackage{nicefrac}
\usepackage[T1]{fontenc}
\usepackage{latexsym,amssymb,amsmath,amsfonts,amsthm}
\usepackage{graphicx}
\usepackage{caption}
\usepackage[utf8]{inputenc}
\usepackage{authblk}
\usepackage{bm}
\numberwithin{equation}{section}
\begin{document}
	
\centerline{\bf\Large A new expression for the Moore--Penrose inverse of a class of matrices}
	
\bigskip

\centerline{Xuefeng Xu}

\medskip

\centerline{\small LSEC, Institute of Computational Mathematics and Scientific/Engineering Computing,}
\centerline{\small Academy of Mathematics and Systems Science, Chinese Academy of Sciences, Beijing 100190, China.}
\centerline{\small  E-mail address: xuxuefeng@lsec.cc.ac.cn}

\medskip

\noindent{\bf \large Abstract}

An expression for the Moore--Penrose inverse of a matrix  of the form $M=XNY$, where $X$ and $Y$ are nonsingular, has been recently established by Castro-Gonz\'{a}lez et al. [1, Theorem~2.2]. The expression plays an essential role in developing explicit expressions for the Moore--Penrose inverse of a two-by-two block matrix. In this paper, we present a new expression for the Moore--Penrose inverse of this class of matrices, which improves the result in [1].

\medskip


\bigskip
\bigskip

\noindent{\bf \large 1.\ Introduction}

\medskip

We first introduce some notations and concepts which are frequently used in the subsequent content. Let $\mathbb{N}^{+}$ and $\mathbb{C}$ denote the set of all positive integers and the field of complex numbers, respectively. Let $\mathbb{C}^{m\times n}$ be the set of all $m\times n$ complex matrices. The identity matrix of order $n$ is denoted by $I_{n}$ or $I$ when its size is clear in the context. For a matrix $A\in\mathbb{C}^{m\times n}$, $A^{\ast}$ denotes the conjugate transpose of $A$. We denote by $\mathcal{R}(A)$ and $\mathcal{N}(A)$ the range and null space of $A$, respectively, namely,  $\mathcal{R}(A):=\{y\in\mathbb{C}^{m}: y=Ax, x\in\mathbb{C}^{n}\}$ and $\mathcal{N}(A):=\{x\in\mathbb{C}^{n}: Ax=0\}$. The \emph{Moore--Penrose inverse} of $A$ is denoted by $A^{\dagger}$, which is defined as the unique matrix $Z\in\mathbb{C}^{n\times m}$ satisfying the following equations:
\begin{displaymath}
(\text{a})\ AZA=A, \quad (\text{b})\ ZAZ=Z, \quad (\text{c})\ (AZ)^{\ast}=AZ,\quad (\text{d})\ (ZA)^{\ast}=ZA.
\end{displaymath} 
The symbols $E_{A}:=I-AA^{\dagger}$ and $F_{A}:=I-A^{\dagger}A$ stand for the orthogonal projectors onto $\mathcal{N}(A^{\ast})$ and $\mathcal{N}(A)$, respectively. A matrix $Z\in\mathbb{C}^{n\times m}$ is called an \emph{inner inverse} of $A$ if it satisfies the equality (a).  

For a matrix $M\in\mathbb{C}^{m\times n}$ which can be decomposed as $M=XNY$, where $X\in\mathbb{C}^{m\times m}$ and $Y\in\mathbb{C}^{n\times n}$ are nonsingular, the equality $M^{\dagger}=Y^{-1}N^{\dagger}X^{-1}$ may fail. Several conditions validating $M^{\dagger}=Y^{-1}N^{\dagger}X^{-1}$ are presented in [2]. Recently, Castro-Gonz\'{a}lez et al. [1] obtained an explicit expression for  $M^{\dagger}$, provided that $XE_{N}=E_{N}$ and $F_{N}Y=F_{N}$. More concretely, it is proved by Castro-Gonz\'{a}lez et al. [1, Theorem 2.2] that 
\begin{align*}
M^{\dagger}=(I+L_{0}^{\ast})(I+L_{0}L_{0}^{\ast})^{-1}Y^{-1}N^{\dagger}X^{-1}(I+R_{0}^{\ast}R_{0})^{-1}(I+R_{0}^{\ast}),\tag{1.1}
\end{align*}
where $R_{0}:=E_{N}(I-X^{-1})$ and $L_{0}:=(I-Y^{-1})F_{N}$. The expression (1.1) is a crucial result in~[1], which can be exploited to establish explicit expressions for the Moore--Penrose inverse of a two-by-two block matrix.

Assume that the singular value decomposition (SVD) of $N\in\mathbb{C}^{m\times n}$ is  $N=U\begin{pmatrix}
\Sigma & 0\\
0 & 0
\end{pmatrix}V^{\ast}$, where $\Sigma\in\mathbb{C}^{r\times r}$ is a diagonal matrix with positive diagonal entries, $r$ is the rank of $N$, and both $U\in\mathbb{C}^{m\times m}$ and $V\in\mathbb{C}^{n\times n}$ are unitary. Let $X\in\mathbb{C}^{m\times m}$ and $Y\in\mathbb{C}^{n\times n}$. We now give two assumptions {\bf A}$_{1}$ and {\bf A}$_{2}$ as follows:
\begin{displaymath}
{\bf A}_{1}: \ X=U\begin{pmatrix}
X_{1} & 0\\
X_{2} & X_{4}
\end{pmatrix}U^{\ast},
\end{displaymath}
where $X_{1}\in\mathbb{C}^{r\times r}$, $X_{2}\in\mathbb{C}^{(m-r)\times r}$, and $X_{4}\in\mathbb{C}^{(m-r)\times (m-r)}$;
\begin{displaymath}
{\bf A}_{2}: \ Y=V\begin{pmatrix}
	Y_{1} & Y_{3}\\
	0 & Y_{4}
\end{pmatrix}V^{\ast},
\end{displaymath}
where $Y_{1}\in\mathbb{C}^{r\times r}$, $Y_{3}\in\mathbb{C}^{r\times (n-r)}$, and $Y_{4}\in\mathbb{C}^{(n-r)\times (n-r)}$. 

In this paper, we further investigate explicit expressions for the Moore--Penrose inverse of this class of matrices. A new expression under weakened conditions for $M^{\dagger}$ is derived, which has enhanced the expression (1.1). More specifically, if the assumptions {\bf A}$_{1}$ and {\bf A}$_{2}$ are satisfied, then we have
\begin{align*}
M^{\dagger}=(I+L^{\ast})(I+LL^{\ast})^{-1}N^{\dagger}N(Y^{-1}N^{\dagger}X^{-1})NN^{\dagger}(I+R^{\ast}R)^{-1}(I+R^{\ast}),\tag{1.2}
\end{align*}
where $R:=XE_{N}X^{-1}(E_{N}-I)$ and $L:=(F_{N}-I)Y^{-1}F_{N}Y$.

The rest of this paper is organized as follows. In Section 2, we first introduce a useful lemma which gives an explicit expression for the Moore--Penrose inverse of a two-by-two block matrix, and then give some specific conditions to validate {\bf A}$_{1}$ and {\bf A}$_{2}$. In Section 3, we present a new and improved expression (i.e., (1.2)) for $M^{\dagger}$ based on the assumptions {\bf A}$_{1}$ and {\bf A}$_{2}$.

\bigskip
\medskip

\noindent{\bf \large 2.\ Preliminaries}

\medskip

In this section, we first introduce a useful lemma, which provides an explicit expression for the Moore--Penrose inverse of a  two-by-two block matrix; see [3]. It is worth mentioning that some improved results of this lemma can be found in [1]. 

\medskip

\noindent{\bf Lemma 2.1.} \emph{Let $M$ be a two-by-two block matrix as the form $M=\begin{pmatrix}
A & C\\
B & D
\end{pmatrix}$. Assume that $\mathcal{R}(B^{\ast})\subseteq\mathcal{R}(A^{\ast})$, $\mathcal{R}(C)\subseteq\mathcal{R}(A)$, and $D-BA^{\dagger}C=0$. Then $M^{\dagger}$ can be given by}
\begin{displaymath}
M^{\dagger}=\begin{pmatrix}
I\\
(A^{\dagger}C)^{\ast}
\end{pmatrix}\Psi A^{\dagger}\Phi \begin{pmatrix}
I & (BA^{\dagger})^{\ast}
\end{pmatrix},
\end{displaymath}
\emph{where $\Phi=\left(I+(BA^{\dagger})^{\ast}BA^{\dagger}\right)^{-1}$ and $\Psi=\left(I+A^{\dagger}C(A^{\dagger}C)^{\ast}\right)^{-1}$.}

\bigskip

Next, we give several specific conditions to guarantee the assumptions {\bf A}$_{1}$ and {\bf A}$_{2}$. 

\medskip

\noindent{\bf Lemma 2.2.} \emph{Let $N\in\mathbb{C}^{m\times n}$ have the singular value decomposition $N=U\begin{pmatrix}
\Sigma & 0\\
0 & 0
\end{pmatrix}V^{\ast}$, where $\Sigma\in\mathbb{C}^{r\times r}$ is a diagonal matrix with positive diagonal entries, $r$ is the rank of $N$, and both $U\in\mathbb{C}^{m\times m}$ and $V\in\mathbb{C}^{n\times n}$ are unitary. Let $X\in\mathbb{C}^{m\times m}$ be an arbitrary matrix. Suppose that one of the following conditions holds:}
\begin{align*}
&{\bf C}_{1}: NN^{\ast}X \ is \ normal;\\
&{\bf C}_{2}: For\ any\ 0\neq c_{1}\in\mathbb{C},\ there \ exists\ k_{1}\in\mathbb{N}^{+}\ such \ that\ (NN^{\ast}X)^{k_{1}}=c_{1}NN^{\dagger};\\
&{\bf C}_{3}: For\ any\ 0\neq c_{2}\in\mathbb{C}\ and\ \ell\in\mathbb{N}^{+},\ there \ exists\ k_{2}\in\mathbb{N}^{+}\ such \ that\ (NN^{\ast}X)^{k_{2}}=c_{2}(NN^{\ast})^{\ell};\\
&{\bf C}_{4}: XE_{N}\ is\ normal;\\
&{\bf C}_{5}: For\ any\ 0\neq c_{3}\in\mathbb{C},\ there \ exists\ k_{3}\in\mathbb{N}^{+}\ such \ that\ (XE_{N})^{k_{3}}=c_{3}E_{N};\\
&{\bf C}_{6}: NN^{\dagger}XE_{N}=0;\\
&{\bf C}_{7}: There\ exists\ k_{4}\in\mathbb{N}^{+}\ such\ that \ (NN^{\ast})^{k_{4}}XE_{N}=0.
\end{align*}
\emph{Then $X$ must be of the form
\begin{displaymath}
X=U\begin{pmatrix}
X_{1} & 0\\
X_{2} & X_{4}
\end{pmatrix}U^{\ast},
\end{displaymath}
where $X_{1}\in\mathbb{C}^{r\times r}$, $X_{2}\in\mathbb{C}^{(m-r)\times r}$, and $X_{4}\in\mathbb{C}^{(m-r)\times (m-r)}$.}

\medskip

\noindent{\bf Proof.} Based on the SVD of $N$, the expressions of $N^{\dagger}$ and $E_{N}$ can be given by
\begin{displaymath}
N^{\dagger}=V\begin{pmatrix}
\Sigma^{-1} & 0\\
0 & 0
\end{pmatrix}
U^{\ast} \quad \text{and} \quad E_{N}=U\begin{pmatrix}
0 & 0\\
0 & I
\end{pmatrix}U^{\ast}.
\end{displaymath}
Partition $U^{\ast}XU$ as $U^{\ast}XU=\begin{pmatrix}
X_{1} & X_{3}\\
X_{2} & X_{4}
\end{pmatrix}$, where $X_{1}\in\mathbb{C}^{r\times r}$, $X_{2}\in\mathbb{C}^{(m-r)\times r}$, $X_{3}\in\mathbb{C}^{r\times (m-r)}$, and $X_{4}\in\mathbb{C}^{(m-r)\times (m-r)}$. Then $X=U\begin{pmatrix}
X_{1} & X_{3}\\
X_{2} & X_{4}
\end{pmatrix}U^{\ast}$.

(i) The condition {\bf C}$_{1}$ states that
\begin{displaymath}
NN^{\ast}X=U\begin{pmatrix}
\Sigma^{2}X_{1} & \Sigma^{2}X_{3}\\
0 & 0
\end{pmatrix}U^{\ast}
\end{displaymath}
is normal, which yields that $\Sigma^{2}X_{1}$ is normal and $\Sigma^{2}X_{3}=0$. It follows from the non-singularity of $\Sigma$ that $X_{3}=0$.

(ii) We have known that $NN^{\ast}X=U\begin{pmatrix}
\Sigma^{2}X_{1} & \Sigma^{2}X_{3}\\
0 & 0
\end{pmatrix}U^{\ast}$. Then, for any $k_{1}\in\mathbb{N}^{+}$, we have  
\begin{displaymath}
(NN^{\ast}X)^{k_{1}}=U\begin{pmatrix}
(\Sigma^{2}X_{1})^{k_{1}} & (\Sigma^{2}X_{1})^{k_{1}-1}\Sigma^{2}X_{3}\\
0 & 0
\end{pmatrix}U^{\ast}.
\end{displaymath}
In addition, it is easy to see that
\begin{align*}
c_{1}NN^{\dagger}=U\begin{pmatrix}c_{1}I & 0\\
0 & 0
\end{pmatrix}U^{\ast}.
\end{align*}
Hence, {\bf C}$_{2}$ implies that $(\Sigma^{2}X_{1})^{k_{1}}=c_{1}I$ and $(\Sigma^{2}X_{1})^{k_{1}-1}\Sigma^{2}X_{3}=0$. Due to the facts that $c_{1}\neq0$ and $\Sigma$ is nonsingular, it follows that $X_{1}$ is nonsingular and $X_{3}=0$.

(iii) Direct calculation yields
\begin{align*}
(NN^{\ast}X)^{k_{2}}&=U\begin{pmatrix}
(\Sigma^{2}X_{1})^{k_{2}} & (\Sigma^{2}X_{1})^{k_{2}-1}\Sigma^{2}X_{3}\\
0 & 0
\end{pmatrix}U^{\ast},\\
c_{2}(NN^{\ast})^{\ell}&=U\begin{pmatrix}
c_{2}\Sigma^{2\ell} & 0\\
0 & 0
\end{pmatrix}U^{\ast}.
\end{align*}
Because $c_{2}\neq0$ and $\Sigma$ is nonsingular, we deduce from {\bf C}$_{3}$ that $\Sigma^{2}X_{1}$ is nonsingular and $(\Sigma^{2}X_{1})^{k_{2}-1}\Sigma^{2}X_{3}=0$. Hence, $X_{3}=0$.

(iv) Straightforward calculation shows
\begin{displaymath}
XE_{N}=U\begin{pmatrix}
0 & X_{3}\\
0 & X_{4}
\end{pmatrix}U^{\ast}.
\end{displaymath}
If $XE_{N}$ is normal, then we get that $X_{4}$ is normal and $X_{3}=0$.

(v) Direct computation yields
\begin{displaymath}
(XE_{N})^{k_{3}}=U\begin{pmatrix}
0 & X_{3}X_{4}^{k_{3}-1}\\
0 & X_{4}^{k_{3}}
\end{pmatrix}U^{\ast}.
\end{displaymath}
It follows from {\bf C}$_{5}$ that $X_{4}^{k_{3}}=c_{3}I$ and $X_{3}X_{4}^{k_{3}-1}=0$. By $c_{3}\neq0$, we derive that $X_{4}$ is nonsingular. Hence, we obtain from $X_{3}X_{4}^{k_{3}-1}=0$ that $X_{3}=0$. 

(vi) It is easy to compute that
\begin{displaymath}
NN^{\dagger}XE_{N}=U\begin{pmatrix}
0 & X_{3}\\
0 & 0
\end{pmatrix}U^{\ast}.
\end{displaymath}
Therefore, $NN^{\dagger}XE_{N}=0$ if and only if $X_{3}=0$.

(vii) Direct calculation yields
\begin{displaymath}
(NN^{\ast})^{k_{4}}XE_{N}=U\begin{pmatrix}
0 & \Sigma^{2k_{4}}X_{3}\\
0 & 0
\end{pmatrix}U^{\ast}.
\end{displaymath}
Due to the fact that $\Sigma$ is nonsingular, it follows that $(NN^{\ast})^{k_{4}}XE_{N}=0$ is equivalent to $X_{3}=0$.

Consequently, if one of the conditions {\bf C}$_{1}$--{\bf C}$_{7}$ holds, then $X$ must be of the from
\begin{displaymath}
X=U\begin{pmatrix}
X_{1} & 0\\
X_{2} & X_{4}
\end{pmatrix}U^{\ast},
\end{displaymath}
which completes the proof. \qed

\bigskip

Analogously, we can prove the following lemma. Its detailed proof is omitted due to limited space.

\medskip

\noindent{\bf Lemma 2.3.} \emph{
Let $Y\in\mathbb{C}^{n\times n}$ and let $N\in\mathbb{C}^{m\times n}$ be the same as in Lemma 2.2. Assume that one of the following conditions holds:}
\begin{align*}
&{\bf C}'_{1}: YN^{\ast}N\ is \ normal;\\
&{\bf C}'_{2}: For\ any\ 0\neq c_{1}'\in\mathbb{C},\ there \ exists\ k_{1}'\in\mathbb{N}^{+}\ such \ that\ (YN^{\ast}N)^{k_{1}'}=c_{1}'N^{\dagger}N;\\
&{\bf C}'_{3}: For\ any\ 0\neq c_{2}'\in\mathbb{C}\ and\ \ell'\in\mathbb{N}^{+},\ there \ exists\ k_{2}'\in\mathbb{N}^{+}\ such \ that\ (YN^{\ast}N)^{k_{2}'}=c_{2}'(N^{\ast}N)^{\ell'};\\
&{\bf C}'_{4}: F_{N}Y\ is\ normal;\\
&{\bf C}'_{5}: For\ any\ 0\neq c_{3}'\in\mathbb{C},\ there \ exists\ k_{3}'\in\mathbb{N}^{+}\ such \ that\ (F_{N}Y)^{k_{3}'}=c_{3}'F_{N};\\
&{\bf C}'_{6}: F_{N}YN^{\dagger}N=0;\\
&{\bf C}'_{7}: There\ exists\ k_{4}'\in\mathbb{N}^{+}\ such\ that \ F_{N}Y(N^{\ast}N)^{k_{4}'}=0.
\end{align*}
\emph{Then $Y$ must be of the form
\begin{displaymath}
Y=V\begin{pmatrix}
Y_{1} & Y_{3}\\
0 & Y_{4}
\end{pmatrix}V^{\ast},
\end{displaymath}
where $Y_{1}\in\mathbb{C}^{r\times r}$, $Y_{3}\in\mathbb{C}^{r\times (n-r)}$, and $Y_{4}\in\mathbb{C}^{(n-r)\times (n-r)}$.}

\bigskip

\noindent{\bf Remark 2.4.} Notice that Lemma 2.2 (resp., Lemma 2.3) does not need the non-singularity of $X$ (resp., $Y$). In addition, the reader can give other conditions to ensure that {\bf A}$_{1}$ and {\bf A}$_{2}$ hold.

\bigskip
\medskip

\noindent{\bf \large 3.\ Main results}

\medskip

In order to prove our main result, we first consider explicit expressions for $(XN)^{\dagger}$ and $(NY)^{\dagger}$. The following theorem provides two applicable formulas for $M_{1}^{\dagger}$ and $M_{2}^{\dagger}$, where $M_{1}=XN$ and $M_{2}=NY$.

\medskip

\noindent{\bf Theorem 3.1.} \emph{Let $N\in\mathbb{C}^{m\times n}$, $X\in\mathbb{C}^{m\times m}$, $Y\in\mathbb{C}^{n\times n}$, $M_{1}=XN$, and $M_{2}=NY$. Suppose that $X$ and $Y$ are nonsingular.} 

\noindent{\emph{(1) If the assumption {\bf A}$_{1}$ holds, then
\begin{align*}
M_{1}^{\dagger}=N^{\dagger}X^{-1}NN^{\dagger}(I+R^{\ast}R)^{-1}(I+R^{\ast}),\tag{3.1}
\end{align*}
where $R=XE_{N}X^{-1}(E_{N}-I)$.}}

\noindent{\emph{(2) If the assumption {\bf A}$_{2}$ holds, then
\begin{displaymath}
M_{2}^{\dagger}=(I+L^{\ast})(I+LL^{\ast})^{-1}N^{\dagger}NY^{-1}N^{\dagger},
\end{displaymath}
where $L=(F_{N}-I)Y^{-1}F_{N}Y$.}}

\medskip

\noindent{\bf Proof.} (1) The assumption {\bf A}$_{1}$ reads $X=U\begin{pmatrix}
X_{1} & 0\\
X_{2} & X_{4}
\end{pmatrix}U^{\ast}$, where $X_{1}\in\mathbb{C}^{r\times r}$ and $r$ is the rank of $N$. It follows from the non-singularity of $X$ that both $X_{1}\in\mathbb{C}^{r\times r}$ and $X_{4}\in\mathbb{C}^{(m-r)\times (m-r)}$ are nonsingular. We define $R:=XE_{N}X^{-1}(E_{N}-I)$. By simple computation, we can get
\begin{displaymath}
R=U\begin{pmatrix}
0 & 0\\
X_{2}X_{1}^{-1} & 0
\end{pmatrix}U^{\ast}=U\begin{pmatrix}
0 & 0\\
G & 0
\end{pmatrix}U^{\ast},
\end{displaymath}
where $G:=X_{2}X_{1}^{-1}$. Because $U$ and $V$ are unitary matrices and
\begin{displaymath}
M_{1}=XN=U\begin{pmatrix}
X_{1}\Sigma & 0\\
X_{2}\Sigma & 0
\end{pmatrix}V^{\ast},
\end{displaymath}
we obtain
\begin{displaymath}
M_{1}^{\dagger}=V\begin{pmatrix}
X_{1}\Sigma & 0\\
X_{2}\Sigma & 0
\end{pmatrix}^{\dagger}U^{\ast}.
\end{displaymath}
Note that $X_{1}\Sigma$ is nonsingular. Using Lemma 2.1, we obtain
\begin{displaymath}
\begin{pmatrix}
X_{1}\Sigma & 0\\
X_{2}\Sigma & 0
\end{pmatrix}^{\dagger}=\begin{pmatrix}
I\\
0
\end{pmatrix} \Sigma^{-1}X_{1}^{-1}(I+G^{\ast}G)^{-1} \begin{pmatrix}
I & G^{\ast}
\end{pmatrix}.
\end{displaymath}
Hence,
\begin{displaymath}
M_{1}^{\dagger}=V\begin{pmatrix}
I\\
0
\end{pmatrix} \Sigma^{-1}X_{1}^{-1}(I+G^{\ast}G)^{-1} \begin{pmatrix}
I & G^{\ast}
\end{pmatrix}U^{\ast}.
\end{displaymath}
Straightforward computation yields
\begin{align*}
N^{\dagger}X^{-1}&=V\begin{pmatrix}
I\\
0
\end{pmatrix} \Sigma^{-1}X_{1}^{-1}\begin{pmatrix}
I & 0
\end{pmatrix}U^{\ast},\\
NN^{\dagger}(I+R^{\ast}R)^{-1}(I+R^{\ast})&=U\begin{pmatrix}
I\\
0
\end{pmatrix}(I+G^{\ast}G)^{-1}\begin{pmatrix}
I & G^{\ast}
\end{pmatrix}U^{\ast}.
\end{align*}
It can be easily seen that $M_{1}^{\dagger}=N^{\dagger}X^{-1}NN^{\dagger}(I+R^{\ast}R)^{-1}(I+R^{\ast})$ holds.

(2) Applying the formula (3.1) to the matrix $Y^{\ast}N^{\ast}$, we obtain
\begin{displaymath}
(M_{2}^{\ast})^{\dagger}=(N^{\ast})^{\dagger}(Y^{\ast})^{-1}N^{\ast}(N^{\ast})^{\dagger}(I+\widehat{R}^{\ast}\widehat{R})^{-1}(I+\widehat{R}^{\ast}),
\end{displaymath}
where
\begin{displaymath} \widehat{R}=Y^{\ast}E_{N^{\ast}}(Y^{\ast})^{-1}(E_{N^{\ast}}-I)=Y^{\ast}(F_{N})^{\ast}(Y^{-1})^{\ast}(F_{N}-I)^{\ast}.
\end{displaymath}
We define $L:=(F_{N}-I)Y^{-1}F_{N}Y$. Then, 
\begin{displaymath}
(M_{2}^{\dagger})^{\ast}=(M_{2}^{\ast})^{\dagger}=(N^{\dagger})^{\ast}(Y^{-1})^{\ast}N^{\ast}(N^{\dagger})^{\ast}(I+LL^{\ast})^{-1}(I+L).
\end{displaymath}
Therefore, we drive that $M_{2}^{\dagger}=(I+L^{\ast})(I+LL^{\ast})^{-1}N^{\dagger}NY^{-1}N^{\dagger}$. \qed

\bigskip

Using Theorem 3.1, we can easily obtain the following expressions for the orthogonal projectors onto $\mathcal{R}(M_{1})$ and $\mathcal{R}(M_{2}^{\ast})$.

\medskip

\noindent{\bf Corollary 3.2.} \emph{Under the same conditions as in Theorem 3.1.} 

\noindent{\emph{(1) If the assumption {\bf A}$_{1}$ is valid, then}}
\begin{align*}
M_{1}M_{1}^{\dagger}=(I+R)NN^{\dagger}(I+R^{\ast}R)^{-1}(I+R^{\ast}).\tag{3.2}
\end{align*}

\noindent{\emph{(2) If the assumption {\bf A}$_{2}$ is valid, then}}
\begin{align*}
M_{2}^{\dagger}M_{2}=(I+L^{\ast})(I+LL^{\ast})^{-1}N^{\dagger}N(I+L).\tag{3.3}
\end{align*}

\medskip

\noindent{\bf Proof.} According to the equality (3.1), it follows that
\begin{align*}
M_{1}M_{1}^{\dagger}=XNN^{\dagger}X^{-1}NN^{\dagger}(I+R^{\ast}R)^{-1}(I+R^{\ast}).\tag{3.4}
\end{align*}
Notice that 
\begin{align*}
(I+R)NN^{\dagger}=\left(I-NN^{\dagger}+XNN^{\dagger}X^{-1}NN^{\dagger}\right)NN^{\dagger}=XNN^{\dagger}X^{-1}NN^{\dagger}.\tag{3.5}
\end{align*}
Inserting (3.5) into (3.4) gives $M_{1}M_{1}^{\dagger}=(I+R)NN^{\dagger}(I+R^{\ast}R)^{-1}(I+R^{\ast})$. Similarly, we can prove the equality (3.3). \qed

\bigskip

Based on the expressions (3.2) and (3.3) for orthogonal projectors $M_{1}M_{1}^{\dagger}$ and $M_{2}^{\dagger}M_{2}$, we can establish the following main result. 

\medskip

\noindent{\bf Theorem 3.3.} \emph{Let $N\in\mathbb{C}^{m\times n}$, $X\in\mathbb{C}^{m\times m}$, $Y\in\mathbb{C}^{n\times n}$, and $M=XNY$. Assume that $X$ and $Y$ are nonsingular. If the assumptions {\bf A}$_{1}$ and {\bf A}$_{2}$ are satisfied, then
\begin{displaymath}
M^{\dagger}=(I+L^{\ast})(I+LL^{\ast})^{-1}N^{\dagger}N(Y^{-1}N^{\dagger}X^{-1})NN^{\dagger}(I+R^{\ast}R)^{-1}(I+R^{\ast}),
\end{displaymath}
where $R=XE_{N}X^{-1}(E_{N}-I)$ and $L=(F_{N}-I)Y^{-1}F_{N}Y$.} 

\medskip

\noindent{\bf Proof.} Note that $Y^{-1}N^{\dagger}X^{-1}$ is an inner inverse of $M$. Then we have
\begin{displaymath}
M^{\dagger}=M^{\dagger}M(Y^{-1}N^{\dagger}X^{-1})MM^{\dagger}.
\end{displaymath}
Let $M_{1}=XN$ and $M_{2}=NY$. We claim that $MM^{\dagger}=M_{1}M_{1}^{\dagger}$ and $M^{\dagger}M=M_{2}^{\dagger}M_{2}$. In fact, it is clear that $MM^{\dagger}$ is the orthogonal projector onto $\mathcal{R}(M)$. Because $Y$ is nonsingular and $M=M_{1}Y$ , it follows that $\mathcal{R}(M)=\mathcal{R}(M_{1})$. Hence, $MM^{\dagger}$ is also an orthogonal projector onto $\mathcal{R}(M_{1})$. Using the uniqueness of orthogonal projectors, we get that $MM^{\dagger}=M_{1}M_{1}^{\dagger}$. Similarly, we can verify that $M^{\dagger}M=M_{2}^{\dagger}M_{2}$. Therefore, we have 
\begin{displaymath}
M^{\dagger}=M_{2}^{\dagger}M_{2}(Y^{-1}N^{\dagger}X^{-1})M_{1}M_{1}^{\dagger}.
\end{displaymath}

Under the assumptions of this theorem, by Corollary 3.2, we have
\begin{displaymath}
M^{\dagger}=(I+L^{\ast})(I+LL^{\ast})^{-1}N^{\dagger}N(I+L)Y^{-1}N^{\dagger}X^{-1}(I+R)NN^{\dagger}(I+R^{\ast}R)^{-1}(I+R^{\ast}).
\end{displaymath}
Using $R=XE_{N}X^{-1}(E_{N}-I)$ and $L=(F_{N}-I)Y^{-1}F_{N}Y$, we obtain
\begin{align*}
(I+L)Y^{-1}N^{\dagger}X^{-1}(I+R)&=Y^{-1}N^{\dagger}X^{-1}+Y^{-1}N^{\dagger}X^{-1}R+LY^{-1}N^{\dagger}X^{-1}+LY^{-1}N^{\dagger}X^{-1}R\\
&=Y^{-1}N^{\dagger}X^{-1},
\end{align*}
where we have applied the facts that $N^{\dagger}E_{N}=0$ and $F_{N}N^{\dagger}=0$. Consequently, we infer that
\begin{displaymath}
M^{\dagger}=(I+L^{\ast})(I+LL^{\ast})^{-1}N^{\dagger}N(Y^{-1}N^{\dagger}X^{-1})NN^{\dagger}(I+R^{\ast}R)^{-1}(I+R^{\ast}).
\end{displaymath}
This completes the proof. \qed

\bigskip

\noindent{\bf Corollary 3.4.} \emph{Under the same conditions as in Theorem 3.3. If both $XE_{N}$ and $F_{N}Y$ are Hermitian, then}
\begin{align*}
M^{\dagger}=(I+L^{\ast})(I+LL^{\ast})^{-1}Y^{-1}N^{\dagger}X^{-1}(I+R^{\ast}R)^{-1}(I+R^{\ast}).\tag{3.6}
\end{align*}

\medskip

\noindent{\bf Proof.} Because $XE_{N}$ and $F_{N}Y$ are Hermitian, by Lemmas 2.2 and 2.3, the assumptions {\bf A}$_{1}$ and {\bf A}$_{2}$ are clearly satisfied. An application of Theorem 3.3 gives
\begin{align*}
M^{\dagger}=(I+L^{\ast})(I+LL^{\ast})^{-1}N^{\dagger}N(Y^{-1}N^{\dagger}X^{-1})NN^{\dagger}(I+R^{\ast}R)^{-1}(I+R^{\ast}).\tag{3.7}
\end{align*}
Due to both $XE_{N}$ and $F_{N}Y$ are Hermitian, it follows that $XE_{N}=E_{N}X^{\ast}$ and $F_{N}Y=Y^{\ast}F_{N}$. Then, $E_{N}(X^{\ast})^{-1}=X^{-1}E_{N}$ and $(Y^{\ast})^{-1}F_{N}=F_{N}Y^{-1}$. Notice that
\begin{displaymath}
N^{\dagger}N(Y^{-1}N^{\dagger}X^{-1})NN^{\dagger}=Y^{-1}N^{\dagger}X^{-1}-F_{N}Y^{-1}N^{\dagger}X^{-1}-Y^{-1}N^{\dagger}X^{-1}E_{N}+F_{N}Y^{-1}N^{\dagger}X^{-1}E_{N}.
\end{displaymath}
Using  $(Y^{\ast})^{-1}F_{N}=F_{N}Y^{-1}$ and $F_{N}N^{\dagger}=0$, we can derive that $F_{N}Y^{-1}N^{\dagger}X^{-1}=0$. By  $E_{N}(X^{\ast})^{-1}=X^{-1}E_{N}$ and $N^{\dagger}E_{N}=0$, we have $Y^{-1}N^{\dagger}X^{-1}E_{N}=0$. Consequently,
\begin{align*}
N^{\dagger}N(Y^{-1}N^{\dagger}X^{-1})NN^{\dagger}=Y^{-1}N^{\dagger}X^{-1}.\tag{3.8}
\end{align*}
By substituting (3.8) into (3.7), we obtain the formula (3.6). \qed

\bigskip

\noindent{\bf Remark 3.5.} If $XE_{N}=E_{N}$ and $F_{N}Y=F_{N}$, the conditions in Corollary 3.4 are obviously satisfied because $E_{N}$ and $F_{N}$ are orthogonal projectors. In this case, \begin{align*}
R&=XE_{N}X^{-1}(E_{N}-I)=E_{N}(X^{-1}E_{N}-X^{-1})=E_{N}(E_{N}-X^{-1})=E_{N}(I-X^{-1})=R_{0},\\
L&=(F_{N}-I)Y^{-1}F_{N}Y=(F_{N}Y^{-1}-Y^{-1})F_{N}=(F_{N}-Y^{-1})F_{N}=(I-Y^{-1})F_{N}=L_{0},
\end{align*}
where $R_{0}$ and $L_{0}$ are defined as in expression (1.1). Therefore, Corollary 3.4 has extended the expression (1.1).

\bigskip
\medskip

\noindent{\bf \large Acknowledgements}

\medskip

The author would like to thank Professor Chen-Song Zhang for his helpful suggestions. This work was supported partially by the National Natural Science Foundation of China Grant 91430215 and 91530323.

\bigskip
\bigskip

\noindent{\bf \large References}

\medskip

\small
{
		
[1]\ N. Castro-Gonz\'{a}lez, M. F. Mart\'{i}nez-Serrano, J. Robles, Expressions for the Moore--Penrose inverse of block matrices involving the Schur complement, Linear Algebra Appl. 471 (2015) 353--368.

[2]\ Y. Tian, The Moore--Penrose inverses of $m\times n$ block matrices and their applications, Linear Algebra Appl. 283 (1998) 35--60.

[3]\ B. Noble, A method for computing the generalized inverse of a matrix, SIAM J. Numer. Anal. 3 (1966) 582--584.

}

\end{document}